\documentclass{amsart}
\usepackage[alphabetic, msc-links, backrefs]{amsrefs}
\usepackage{amssymb}
\usepackage{amsthm}
\usepackage{apptools}
\usepackage{chngcntr}
\usepackage{color}
\usepackage{enumerate}
\usepackage{graphicx}
\usepackage{hyperref}
\usepackage{verbatim}
\usepackage[dvipsnames]{xcolor}

\newtheorem{theorem}            {Theorem}       

\newtheorem{corollary}          [theorem]       {Corollary}

\newtheorem*{theorem*}                          {Theorem}
\theoremstyle{definition}
\newtheorem{definition}         [theorem]       {Definition}

\theoremstyle{definition}
\newtheorem{remark}             [theorem]       {Remark}

\allowdisplaybreaks[4]   

\title[Type II Singularities of LMCF with Zero Maslov Class]{Type II Singularities of Lagrangian Mean Curvature Flow with Zero Maslov Class}

\author{Xiang Li$^1$}
\author{Yong Luo$^2$}
\author{Jun Sun$^{3}$}

\address{$^1$School of Mathematical Sciences, Jiangsu University, Zhenjiang 212013, China}
\email{lihsiang@mail.ustc.edu.cn}
\address{$^2$Mathematical Science Research Center, Chongqing University of Technology, Chongqing 400054, China}
\email{yongluo-math@cqut.edu.cn}
\address{$^3$School of Mathematics and Statistics, Wuhan University, Wuhan 430072, China}
\email{sunjun@whu.edu.cn}

\begin{document}
	
\maketitle

\begin{abstract}
In this paper, we will prove some rigidity theorems for blow up limits to Type II singularities of Lagrangian mean curvature flow with zero Maslov class or almost calibrated Lagrangian mean curvature flows, especially for Lagrangian translating solitons in arbitrary dimension. These theorems generalized previous corresponding results from two dimensional case to arbitrarily dimensional case. 

\vspace{.1in}

\end{abstract}

\section{Introduction}

\vspace{.1in}

Let $M^n$ be a Calabi-Yau manifold of complex dimension $n$ with a K\"ahler form $\omega$, a complex structure $J$, a K\"ahler metric $g$ and a parallel holomorphic $(n,0)$-form $\Omega$ of unit length.  

An immersed submanifold $\Sigma$ in $M$ is Lagrangian if $\omega|_{\Sigma}=0$. 
The induced volume form $d\mu_{\Sigma}$ on a Lagrangian submanifold $\Sigma$ from the Ricci-flat metric $g$ is related to $\Omega$ by (see \cite{HarLaw})
$$\Omega |_\Sigma=e^{i\theta}d\mu_\Sigma,
$$
where $\theta$ is a  multivalued function called the Lagrangian angle.
When the Lagrangian angle is a single valued function the Lagrangian is called {\it zero Maslov class} and if $\cos\theta\geq \delta>0$ for some positive $\delta$, then $\Sigma$ is called {\it almost calibrated}. 
If $\theta\equiv constant$, then $\Sigma$ is a {\it Special Lagrangian}. 

\vspace{.1in}

The existence of Special Lagrangians in a Calabi-Yau manifold is an important problem in differential geometry and mirror symmetry. The Special Lagrangians not only are minimal submanifolds, but also are calibrated submanifolds of Calabi-Yau manifolds (\cite{HarLaw}). Therefore, the Special Lagrangians are area-minimizing in their homology class. Schoen and Wolfson (\cite{ScW}) studied the minimization problem and showed that, when the real dimension is four, a Lagrangian minimizing area among all Lagrangians in a given class exists, is smooth everywhere except finitely many points, but is not necessarily a minimal surface. Later Wolfson (\cite{Wolf}) found a Lagrangian sphere with nontrivial homology on a given $\mathbf{K}3$ surface such that the Lagrangian that minimizes area among all Lagrangians in this class is not a Special Lagrangian and the surface that minimizes area among all surfaces in this class is not Lagrangian. This shows the subtle nature of the problem (\cite{NevesAndre2}). In addition to using variational methods to find minimal surfaces, another important idea is to use the mean curvature flow, which is the negative gradient flow for the area functional.

Smoczyk (\cite{Smo}) proved that the ``Lagrangian" condition is preserved by mean curvature flow, in which case the flow is called the Lagrangian mean curvature flow.
It also follows from the maximum principle that the ``zero Maslov class" property or ``almost calibrated" condition is preserved by the Lagrangian mean curvature flow, in which case the flow is called Lagrangian mean curvature flow with zero Maslov class, or almost calibrated Lagrangian mean curvature flow.

By estalbishing a new monotonicity formula, Chen-Li (\cite{CL2}) and Wang (\cite{Wang1}) independently proved that there is no finite time Type I singularity for almost calibrated Lagrangian mean curvature flow. Neves (\cite{NevesAndre}) proved the same conclusion for Lagrangian mean curvature flow with zero Maslov class. Furthermore, Neves constructed examples of Lagrangian mean curvature flow with zero Maslov class which develops finite time singularity (\cite{NevesAndre}, \cite{NevesAndre2}). Therefore, it is important to study blow up flows for Type II singularities of the Lagrangian mean curvature flow, which is a nonflat complete eternal solution  (i.e. a solution defined on the whole time interval $(-\infty, +\infty)$) to the Lagrangian mean curvature in ${\mathbb C}^n$ with bounded second fundamental form, at most Euclidean area growth, and bounded Lagrangian angle (or almost calibrated).

\vspace{.1in}

In $\mathbb{C}^{2}$, Han-Li-Sun (\cite{HLS}) have studied the properties of the general limit flow $\Sigma^{\infty}_{s}$ of the almost calibrated Lagrangian mean curvature flow. In particular, they provided an estimate for the norm of the mean curvature vector in terms of the oscillation of the Lagrangian angle for the limit flow. The first result in this paper is to generalize this result for Lagrangian mean curvature flow with zero Maslov class in any dimension $n$.
Specifically, we have:

\begin{theorem}\label{thm1.1}
	Suppose that $\Sigma_{t}$ $(t\in(-\infty,0])$ is a complete proper Lagrangian mean curvature flow in $\mathbb{C}^{n}$ with bounded Lagrangian angle, i.e., $\sup_{t\in (-\infty,0]}\sup_{\Sigma_t}|\theta|\leq \Lambda$ for some positive constant $\Lambda$. Assume further that $\sup_{t\in(-\infty,0]}\sup_{\Sigma_{t}}|A|^{2}=1$. Then we have
	\begin{equation*}
		h^{2}=\sup_{t\in(-\infty,0]}\sup_{\Sigma_{t}}|H|^{2}\leq\big(\sup_{t\in(-\infty,0]}\sup_{\Sigma_{t}}\theta-\inf_{t\in(-\infty,0]}\inf_{\Sigma_{t}}\theta\big)^{2}.
	\end{equation*}
Here and in the sequel, $H$ and $A$ denote the mean curvature vector and the second fundamental form of $\Sigma_t$ in ${\mathbb C}^n$, respectively.
\end{theorem}

\vspace{.1in}

In the same paper, Han-Li-Sun (\cite{HLS}) also proved that an eternal solution to the almost calbrated Lagrangian mean curvature flow in ${\mathbb C}^2$ which is flat at all time must be flat planes, hence cannot arise as blow up limit of almost calibrated Lagrangian mean curvature flow. By the Gauss equation, we can see that a surface is flat (i.e. the Gauss curvature vanishes everywhere) if and only if $|H|^2=|A|^2$ on $\Sigma$. Recently, Li-Sun (\cite{LS}) have shown that a stronger result that any almost calibrated Lagrangian eternal mean curvature flow in $\mathbb{C}^{2}$ with $\cos\theta\geq\delta>0$ and $|H|^{2}\geq\varepsilon|A|^{2}$ with $\varepsilon>1-\delta$ must be flat planes.
Our next theorem generalizes it to arbitrary dimension $n$.

\begin{theorem}\label{thm1.2}
	Any proper almost calibrated Lagrangian eternal mean curvature flow in $\mathbb{C}^{n}$ with $\cos\theta\geq\delta>0$ and $|H|^{2}\geq\varepsilon|A|^{2}$ with $\varepsilon>1-\delta$ must be flat planes.
\end{theorem}

\vspace{.1in}

As a consequence, we also have the following improvement in (\cite{HLS}), (\cite{HanS}) and (\cite{LS}).

\begin{corollary}\label{cor2}
	Any proper eternal mean curvature flow in ${\mathbb C}^{n}$ with scalar curvature $R_{\Sigma_t}\geq-(1-\varepsilon)|A|^{2}$ for $\varepsilon>1-\delta$ cannot arise as blow up flow of almost calibrated Lagrangian mean curvature flow.
\end{corollary}

\vspace{.1in}

An important type of eternal solutions to the Lagrangian mean curvature flow are Lagrangian translating solitons.  
There are some nonexistence results and rigidity theorems on Lagrangian translating solitons (e.g., \cite{HL}, \cite{HanS}, \cite{NT}, \cite{QiuHongbing}, \cite{Sun1,Sun2,Sun3}, etc.). 

Recall that a surface $\Sigma^{n}$ in $\mathbb{R}^{n+k}$ is called a translating soliton (or translator) of the mean curvature flow, if it satisfies
\begin{equation*}
	T^{\bot}=H,
\end{equation*}
where $H$ is the mean curvature vector of $\Sigma^{n}$ in $\mathbb{R}^{n+k}$ and $T$ is an unit constant vector in $\mathbb{R}^{n+k}$. Let $V$ be the tangent part of $T$. 
Then we have
\begin{equation}\label{e-TS}
	T=V+H.
\end{equation}

In \cite{Sun2}, Sun proved that the infimum of $|H|$ on any complete almost-calibrated Lagrangian translating soliton in ${\mathbb C}^2$ with $\cos\theta\geq\delta>0$ and quadratic area growth must be zero.
Our next result is to improve this result for complete Lagrangian translating soliton with zero Maslov class of any dimension:

\begin{theorem}\label{thm1.3}
	Suppose $\Sigma$ is a complete Lagrangian translating soliton with zero Maslov class and has Euclidean area growth in $\mathbb{C}^{n}$. Then:
	\begin{equation*}
		\inf_{\Sigma}|H|^{2}=0.
	\end{equation*}
\end{theorem}

\vspace{.1in}

It is well known that a translating soliton can be viewed as a critical point of the
functional
\begin{equation*}
	\mathcal{L}(\Sigma)=\int_{\Sigma}e^{\left\langle T,{\bf{x}}\right\rangle}d\mu,
\end{equation*}
where ${\bf{x}}$ is the position vector in $\mathbb{R}^{n+k}$, and $d\mu$ is the volume form on $\Sigma$ induced from the Euclidean space $\mathbb{R}^{n+k}$. 
For our convenience, we denote $d\widetilde{\mu}=e^{\left\langle T,{\bf{x}}\right\rangle}d\mu$.

It is easy to see that a complete translating soliton cannot be compact.
Inspired by Xin (\cite{Xin}) and Sun's (\cite{Sun3}) work, we have

\begin{theorem}\label{thm1.4}
	Suppose $\Sigma$ is a complete Lagrangian translating soliton with zero Maslov class in $\mathbb{C}^{n}$ and mean curvature vector $H\in L^{1}(d\widetilde{\mu})$. Then $\Sigma$ must be a minimal Lagrangian submanifold.
\end{theorem}

\begin{remark}
	In Theorem \ref{thm1.4} and Theorem \ref{thm1.5} below, we can only prove that the translating solitons with suitable restrictions are  minimal instead of flat plane. Actually, when $n=2$, it is easy to see that a minimal Lagrangian translating soliton in ${\mathbb C}^2$ must be a plane (see (\ref{AHV})). However, when $n\geq 3$, a minimal Lagrangian translating soliton in ${\mathbb C}^n$ need not be a flat plane (\cite{NT}). For instance, if $\Sigma$ is any nonflat Special Lagrangian in ${\mathbb C}^2$, then $\tilde{\Sigma}:={\mathbb R}\times \Sigma\subset {\mathbb C}^3$ is minimal Lagrangian translating soliton but nonflat.
\end{remark}

\begin{corollary}\label{cor2}
    Any translating soliton with mean curvature vector $H\in L^{1}(d\widetilde{\mu})$ in $\mathbb{C}^{2}$  cannot arise as blow up limit of the Lagrangian mean curvature flow with zero Maslov class.
\end{corollary}

In order to state the next theorem, we give the definition of weighted polynomial area growth:
\begin{definition}
	We say a submanifold $\Sigma^n$ in ${\mathbb R}^{n+k}$ has weighted polynomial area growth, if there is a constant $D_{0}>0$ and $d>0$, such that
	\begin{equation}
		\widetilde{\mu}(\Sigma\cap B(r)):=\int_{\Sigma\cap B(r)}e^{\left\langle T,{\bf{x}}\right\rangle}d\mu
		\leq D_{0}r^{d},
	\end{equation}
	for any $r\geq1$ holds, where $B(r)$ is the ball of radius $r$ in ${\mathbb R}^{n+k}$.
\end{definition}
\vspace{.1in}
Then we have
\begin{theorem}\label{thm1.5}
	Suppose $\Sigma$ is a complete Lagrangian translating soliton with zero Maslov class in $\mathbb{C}^{n}$ and weighted polynomial area growth. Then $\Sigma$ must be a minimal submanifold.
\end{theorem}

\vspace{.1in}

As a corollary, we get that

\begin{corollary}\label{cor1}
   Any translating soliton with weighted polynomial area growth in $\mathbb{C}^{2}$ cannot arise as blow up limit of the Lagrangian mean curvature flow with zero Maslov class.
\end{corollary}

\vspace{.1in}

By the monotonicity formula to the mean curvature flow, we know that the blowup limit of mean curvature flow must have polynomial area growth. However it is not clear whether the blow up limit has weighted polynomial area growth.

\vspace{.1in}

In the end of this section let us emphasize several significance of results and innovations of methodology in this paper. Firstly, all the results proved in this paper are true for any dimension instead of dimension $n=2$. Actually, the dimension 2 is a crucial assumption in \cite{HLS} and \cite{LS}. For $n=2$, they can use the weighted monotonicity formula to estimate the maximum of the mean curvature vector. These arguments are invalid for $n>2$. To overcome this difficulty, instead of deriving the estimate for the mean curvature vector, we choose a new cutoff function to obtain the desired decay. Secondly, all the previous results are proved for almost calibrated Lagrangian translating solitons, where they used $\frac{1}{\cos\theta}$ to construct the test functions. In this paper, we choose to look for cutoff funtions depending on $\theta$ itself so that most of our results are true for Lagrangian translating solitons with zero Maslov class.

\vspace{.1in}

The subsequent sections are organized as follows: in Section 2, we provide some preliminary materials and in Section 3 and Section 4, we prove the main theorems.

\vspace{.2in}

\section{Notations and Preliminaries}

\vspace{.1in}

Let $F_{0}:N^n\rightarrow M^{n+k}$ be a smooth immersion of a smooth manifold $N$. 
The mean curvature flow with initial condition $F_0$ is a smooth family of immersions $F:N^n\times[0,T)\rightarrow M^{n+k}$ satisfying
\begin{eqnarray*}
	\begin{cases}
		\frac{\partial}{\partial t}F(x,t)=H(x.t),\\
		F(x,0)=F_{0}(x),
	\end{cases}
\end{eqnarray*}
where $H(x,t)=H^{\alpha}\nu_{\alpha}=h^{\alpha}_{ii}\nu_{\alpha}$ is the mean curvature vector of the submanifold $\Sigma_{t}=F(N, t)$ at $x$. The mean curvature fow is the negative gradient flow of the area functional. If the flow exists globally and converges at infinity, then the limit must be a minimal submanifold.

Denote by $A$ the second fundamental form of $\Sigma_t$ in $M$, and the Riemannian metric on $M$ by $\langle\cdot,\cdot\rangle$. In normal coordinates around a point in $\Sigma_t$, the induced metric on $\Sigma_t$ from
$\langle\cdot,\cdot\rangle$ is given by
$$
g_{ij}=\langle \partial_iF,\partial_jF\rangle,
$$
where $\partial_i$ ($i=1,\cdots ,n$) are the partial derivatives with
respect to the local coordinates. In the sequel,  we denote by
$\Delta$ the Laplace operator and $\nabla$ the covariant connection associated with the induced metric on $\Sigma_t$ respectively. We choose the orthonormal frame $e_1,\cdots,e_n$, $\nu_1,\cdots,\nu_{m-n}$
of $M$ along $\Sigma_t$ such that $e_1,...,e_n$ are tangent
vectors of $\Sigma_t$ and $v_1,...,v_{m-n}$ are in the normal bundle
over $\Sigma_t$. We can write:
\begin{eqnarray*}
A&=&A^\alpha v_\alpha\\
H&=&-H^\alpha v_\alpha.
\end{eqnarray*}
Let $A^\alpha=(h^\alpha_{ij})$ where $(h^\alpha_{ij})$ is a matrix, by the
Weingarten equation, we have
$$
h_{ij}^\alpha=\langle\bar{\nabla}_iv_\alpha,e_j\rangle
=\langle\bar{\nabla}_jv_\alpha,e_i\rangle=h^\alpha_{ji},
$$
where $\bar{\nabla}$ is the Levi-Civita connection on $M$. The trace
and the norm of the second fundamental form of $\Sigma_t$ in $M$
are:
$$
H^\alpha = g^{ij}h^\alpha_{ij}=h_{ii}^\alpha $$ $$ |A|^2=\sum_{\alpha}
|A^\alpha |^2=g^{ij}g^{kl}h^\alpha_{ik}h^\alpha_{jl}=h_{ik}^\alpha h_{ik}^\alpha.
$$

For a Lagrangian submanifold $N$, we have the following relation between mean curvature vector and the Lagrangian angle (\cite{TY}):
\begin{equation*}\label{HJtheta}
	H=J\nabla\theta.
\end{equation*}
It then follows that
\begin{equation}\label{nablatheta=H}
	|\nabla\theta|=|H|.
\end{equation}

Smoczyk (\cite{Smo}) proved that if the initial submanifold is Lagrangian, then along the mean curvature flow, it remain Lagrangian for positive time. In this case, the flow is called ``Lagrangian mean curvature flow". Smoczyk (\cite{Smo2}) also derived the following equation for the Lagrangian angle along the mean curvature flow:
\begin{equation}\label{theta}
	(\frac{\partial}{\partial t}-\Delta)\theta=0
\end{equation}
and thus
\begin{equation}\label{Htheta}
	(\frac{\partial}{\partial t}-\Delta)\cos\theta =|H|^2\cos\theta.
\end{equation}
As a result, if the initial submanifold is almost calibrated, i.e., $\cos\theta(\cdot,0)$ has a positive lower bound, then by applying the parabolic maximum principle to this evolution equation, one concludes that $\cos\theta$ remains positive as long as the mean curvature flow has a smooth solution, i.e., almost calibrated condition is preserved by mean curvature flow, in which case the flow is called ``almost calibrated Lagrangian mean curvature flow". (\ref{theta}) also implies that ``zero Maslov class" property is preserved along the Lagrangian mean curvature flow. In this case the flow is called ``Lagrangian mean curvature flow with zero Maslov class".

From (\ref{theta}), it is easy to see that on Lagrangian translating solitons, we have:
\begin{equation}\label{Deltatheta}
 	\Delta\theta+\left\langle V,\nabla\theta\right\rangle=0,
\end{equation}
where $V$ is given by (\ref{e-TS}).
This implies that
\begin{equation}\label{Deltacostheta}
	-\Delta\cos\theta=|H|^{2}\cos\theta+V\cdot\nabla\cos\theta,
\end{equation}
which will be used later.

For a two dimensional translating soliton in $\mathbb{R}^{4}$, Han-Li (\cite{HL}) computed the following identities on translating solitons at the points where $V\neq 0$:
\begin{equation}\label{AHV}
	|A|^{2}=|H|^{2}+2\frac{|\nabla H|^{2}}{|V|^{2}}
	+\frac{V\cdot\nabla| H|^{2}}{|V|^{2}}.
\end{equation}
This implies that two dimensional translating soliton in $\mathbb{R}^{4}$ is a flat plane if and only if it is minimal.

Finally, we also have the following well-known evolution equation for the mean curvature vector along mean curvature flow in ${\mathbb R}^{n+k}$ which will be used later (see, for example, (\cite{AB})):
\begin{equation}\label{H2}
	(\frac{\partial}{\partial t}-\Delta)|H|^{2}=-2|\nabla^{\bot}H|^{2}+2\sum_{i,j,\alpha,\beta}H^{\alpha}H^{\beta}h^{\alpha}_{ij}h^{\beta}_{ij}\leq -2|\nabla^{\bot}H|^{2}+2|A|^2|H|^2.
\end{equation}

\vspace{.2in}

\section{Rigidity theorems for blow up flows of Lagrangian mean curvature flows}

\vspace{.1in}

In this section, we will prove the first two main theorems for blow up flows of Lagrangian mean curvature flows.

\vspace{.1in}

\noindent\textbf{Proof of Theorem \ref{thm1.1}:}
Without loss of generality, we may assume
\begin{equation*}
	\inf_{t\in(-\infty,0]}\inf_{\Sigma_{t}}\theta=0.
\end{equation*}
If $h=0$ or $\eta:=\sup_{t\in(-\infty,0]}\sup_{\Sigma_{t}}\theta=0$, it is evident that the result holds. Now we assume that $h>0$, $\eta>0$ and fix any positive number $T$.
 
We introduce a new function $f(X,t)=|H|^{2}+p\theta^{2}$, where $X$ is the position vector in ${\mathbb R}^{2n}$, $p>1$ and $t\in[-T,0]$. 
Using the evolution equation \eqref{theta} for $\theta$ and the evolution equation \eqref{H2} for $|H|^{2}$, we get 
\begin{equation*}
	\begin{split}
			(\Delta-\frac{\partial}{\partial t})f&\geq2|\nabla^{\bot}H|^{2}-2|H|^{2}|A|^{2}+2p|\nabla\theta|^{2}\\
			&\geq2(p-1)|H|^{2}.
	\end{split}
\end{equation*}
Here, we have used the assumption $\sup_{t\in(-\infty,0]}\sup_{\Sigma_{t}}|A|^{2}=1$
and the fact \eqref{nablatheta=H}.

Let $\psi(x)$ be a $C^{2}$ function on $[0,\infty)$ such that
\begin{eqnarray*}
	\psi(x)=\left\{\begin{array}{clcr}  1, &{\rm if} & 0\leq x\leq\frac{1}{2},\\
		0, &{\rm if}  & x\geq 1,
	\end{array}\right.
\end{eqnarray*}
$$ 0\leq\psi(x)\leq 1,~~~ \psi'(x)\leq 0,~~~ |\psi''(x)|\leq C ~~~{\rm and
}~~~ \frac{|\psi'(x)|^2}{\psi(x)}\leq C,
$$ 
where $C$ is an absolute constant.

We construct a cutoff function $g(X,t)=\psi(\frac{|X|^{2}}{R^{2}})\psi(\frac{t}{T})$.
A straightforward computation shows that
\begin{equation*}
	\begin{split}
		(\Delta-\frac{\partial}{\partial t})g&=4\psi^{''}(\frac{|X|^{2}}{R^{2}})\psi(\frac{t}{T})\frac{\left\langle X,\nabla X\right\rangle^{2}}{R^{4}}
		+\psi(\frac{t}{T})\frac{2n\psi^{'}(\frac{|X|^{2}}{R^{2}})}{R^{2}}
		-\frac{1}{T}\psi^{'}(\frac{t}{T})\psi(\frac{|X|^{2}}{R^{2}})\\
		&\geq-\frac{C}{R^{2}}-\frac{C}{T}
	\end{split}
\end{equation*}
and
\begin{equation*}\label{Hgnablag}
	\frac{|\nabla g|^2}{g}\leq\frac{C}{R^2}.
\end{equation*}

Set $\phi=gf$. Since  $\phi\big|_{\partial(\overline{B}_{R}\times[-T,0])}=0$, then $\phi$ must achieve its maximum at some points $(X_{i},t_{i})\in B_{R}\times(-T,0)$.
By the maximum value principle, at $(X_{i},t_{i})$, we obtain that
\begin{equation*}
	\frac{\partial}{\partial t}\phi(X_{i},t_{i})\geq0,
\end{equation*}
\begin{equation*}
	\nabla\phi(X_{i},t_{i})=0,
\end{equation*}
\begin{equation*}
	\Delta\phi(X_{i},t_{i})\leq0,
\end{equation*}
i.e.,
\begin{equation*}
	\big(\Delta-\frac{\partial}{\partial t}\big)\phi=\big(\Delta-\frac{\partial}{\partial t}\big)gf\leq0
\end{equation*}
and
\begin{equation*}
	\nabla g=-\frac{g}{f}\nabla f.
\end{equation*}
Hence, we get
\begin{equation*}
	\begin{split}
		0&\geq\big(\Delta-\frac{\partial}{\partial t}\big)gf\\
		&=f\big(\Delta-\frac{\partial}{\partial t}\big)g+g\big(\Delta-\frac{\partial}{\partial t}\big)f+2\nabla g\nabla f\\
		&\geq -\frac{C}{R^{2}}f-\frac{C}{T}f+2(p-1)g|H|^{2}-\frac{|\nabla g|^{2}}{g}f\\
		&\geq -\frac{C}{R^{2}}f-\frac{C}{T}f+2(p-1)g|H|^{2}.
	\end{split}
\end{equation*}
Since $p>1$ and $f$ is bounded, we get 
\begin{equation*}
	g|H|^{2}(X_{i},t_{i})\leq\frac{C}{(p-1)R^{2}}+\frac{C}{(p-1)T}.
\end{equation*}
Therefore,
\begin{equation*}
	\sup_{B_{\frac{R}{2}}\times [-T,0]}f(X_{i},t_{i})\leq\frac{C}{(p-1)R^{2}}+\frac{C}{(p-1)T}+p\sup_{B_{R}\times [-T,0]}\theta^{2}.
\end{equation*}
Letting $R\to+\infty$ and $T\to+\infty$, we obtain
\begin{equation*}
	h^{2}\leq p\eta^{2}.
\end{equation*}
Letting $p\to 1$, we get the desired inequality.
\hfill $\boxed{}$

\vspace{.1in}

\noindent\textbf{Proof of Theorem \ref{thm1.2}:}
We first prove that the eternal flow $\Sigma_t, t\in (-\infty, +\infty)$ is minimal under the assumption that
\begin{equation}\label{varepsilon}
	|H|^{2}\geq\varepsilon|A|^{2},
\end{equation}
where $\varepsilon>1-\delta$. We prove it by contradiction.

As explained earlier, we suppose on the contrary that the eternal solution is not minimal. Without loss of generality, we may assume that $|H|(0,0)=1$. Let $f$ be a positive increasing function which is to be determined later. We compute $\big(\Delta-\frac{\partial}{\partial t}\big) (|H|^2f(\frac{1}{\cos\theta}))$.
First note that by (\ref{Htheta}),
\begin{equation}\label{e-theta}
	\begin{split}
		\big(\Delta-\frac{\partial}{\partial t}\big) \frac{1}{\cos\theta}&=
		-\frac{\big(\Delta-\frac{\partial}{\partial t}\big) \cos\theta}{\cos^{2}\theta}+2\cos\theta|\nabla\frac{1}{\cos\theta}|^{2} \\
		&=\frac{|H|^{2}}{\cos\theta}+2\cos\theta|\nabla\frac{1}{\cos\theta}|^{2}.
	\end{split}
\end{equation}
On the other hand, the  equality (\ref{H2}) combining with
Cauchy-Schwarz inequality yields
\begin{equation}
	\begin{split}
		\big(\Delta-\frac{\partial}{\partial t}\big) |H|^{2}&=2|\nabla^{\bot}
		H|^{2}-2\sum_{i,j,\alpha,\beta}H^{\alpha}H^{\beta}h^{\alpha}_{ij}h^{\beta}_{ij} \\
		&\geq2|\nabla^{\bot}
		H|^{2}-2|A|^{2}|H|^{2}.
	\end{split}
\end{equation}
We calculate using (\ref{H2}), (\ref{e-theta}) and rearranging the terms to obtain
\begin{align}\label{HH2f}
	\big(\Delta-\frac{\partial}{\partial t}\big)(|H|^2f(\frac{1}{\cos\theta}))
	&=\big(\Delta-\frac{\partial}{\partial t}\big)|H|^2f(\frac{1}{\cos\theta})+|H|^2\big(\Delta-\frac{\partial}{\partial t}\big)(f(\frac{1}{\cos\theta})) \notag \\
	&\quad+2\nabla
	|H|^2\cdot\nabla f(\frac{1}{\cos\theta}) \notag \\
	&\geq  f(\frac{1}{\cos\theta})(2|\nabla^{\bot} H|^{2}-2|A|^{2}|H|^{2}) \notag \\
	& \quad +|H|^{2}\big(f'\frac{|H|^{2}}{\cos\theta}+2f'\cos\theta|\nabla\frac{1}{\cos\theta}|^{2}+f''|\nabla\frac{1}{\cos\theta}|^{2}\big)\\
	&\quad +2\frac{\nabla(f|H|^{2})-|H|^{2}\nabla f}{f}\cdot \nabla f(\frac{1}{\cos\theta}) \notag \\
	&= |H|^{2}f\big(2\frac{|\nabla^{\bot}H|^{2}}{|H|^{2}}-2|A|^{2}
	+\frac{f'}{f}\frac{|H|^{2}}{\cos\theta}\big) \notag \\
	&\quad+|H|^{2}\big(f''-2\frac{(f')^{2}}{f}+2f'\cos\theta\big) |\nabla\frac{1}{\cos\theta}|^{2} \notag \\
	&\quad +2|H|^{2}\frac{\nabla(f|H|^{2})}{f|H|^{2}}\cdot \nabla
	f(\frac{1}{\cos\theta}).  \notag
\end{align}
Set $\phi=f|H|^{2}$. At the point where $\phi \neq 0$, it is easy to
see that
\begin{equation*}
	\nabla \phi= f\nabla|H|^{2} + |H|^{2}f' \nabla\frac{1}{\cos\theta},
\end{equation*}
i.e.,
\begin{equation}\label{nablacosalpha}
	\nabla\frac{1}{\cos\theta}=\frac{f}{f'}\big(\frac{\nabla \phi}{\phi}-\frac{\nabla
		|H|^{2}}{|H|^{2}}\big).
\end{equation}
Plugging (\ref{nablacosalpha}) into (\ref{HH2f}), we obtain
\begin{align}\label{Hphi}
	\big(\Delta-\frac{\partial}{\partial t}\big)\phi
	&\geq \phi\big(2\frac{|\nabla^{\bot} H|^{2}}{|H|^{2}}-2|A|^{2}
	+\frac{f'}{f}\frac{|H|^{2}}{\cos\theta}\big)\notag \\
	&\quad +\frac{\phi f}{(f')^{2}}\big(f''-2\frac{(f')^{2}}{f}+2f'\cos\theta\big)
	\big(\frac{|\nabla |H|^{2}|^{2}}{|H|^{4}}-2\frac{\nabla |H|^{2}}{|H|^{2}}\cdot\frac{\nabla \phi}{\phi}
	+\frac{|\nabla \phi|^{2}}{\phi^{2}}\big)\notag  \\
	&\quad +2|H|^{2}f'\frac{\nabla\phi}{\phi}\nabla\frac{1}{\cos\theta} \notag \\
	&=\phi\big(\frac{f'}{f}\frac{|H|^{2}}{\cos\theta}-2|A|^{2}\big)
	+ \phi\big(2\frac{|\nabla^{\bot} H|^{2}}{|H|^{2}}
	+4\frac{ff''}{(f')^{2}}\frac{|\nabla |H||^{2}}{|H|^{2}}-8\frac{|\nabla |H||^{2}}{|H|^{2}}\notag \\
	&\quad +8\frac{f}{f'}\cos\theta \frac{|\nabla |H||^2}{|H|^{2}}\big)
	+\phi\big(\frac{ff''}{(f')^{2}}+2\frac{f}{f'}\cos\theta-2\big)\big(\frac{|\nabla \phi|^{2}}{\phi^{2}}
	-2\frac{\nabla|H|^{2}}{|H|^{2}}\cdot\frac{\nabla \phi}{\phi}\big)\notag  \\
	&\quad +2|H|^{2}f'\frac{\nabla\phi}{\phi}\nabla \frac{1}{\cos\theta} \\
	&\geq \phi\big(\frac{f'}{f}\frac{|H|^{2}}{\cos\theta}-2|A|^{2}\big)
	+\phi\big(4\frac{ff''}{(f')^{2}}+8\frac{f}{f'}\cos\theta-6\big)\frac{|\nabla |H||^{2}}{|H|^{2}}\notag \\
	&\quad +\phi\big(\frac{ff''}{(f')^{2}}+2\frac{f}{f'}\cos\theta-2\big)\big(\frac{|\nabla \phi|^{2}}{\phi^{2}}
	-2\frac{\nabla|H|^{2}}{|H|^{2}}\cdot\frac{\nabla \phi}{\phi}\big)\notag  \\
	&\quad +2|H|^{2}f'\frac{\nabla\phi}{\phi}\nabla \frac{1}{\cos\theta}\notag  \\
	&= \phi\big(\frac{f'}{f}\frac{|H|^{2}}{\cos\theta}-2|A|^{2}\big)
	+\phi\big(4\frac{ff''}{(f')^{2}}+8\frac{f}{f'}\cos\theta-6\big)\frac{|\nabla |H||^{2}}{|H|^{2}}\notag  \\
	&\quad -\phi\big(\frac{ff''}{(f')^{2}}+2\frac{f}{f'}\cos\theta-2\big)\big(\frac{|\nabla \phi|^{2}}{\phi^{2}}
	-2\frac{f'}{f}\nabla\frac{1}{\cos\theta}\cdot\frac{\nabla \phi}{\phi}\big)\notag  \\
	&\quad +2|H|^{2}f'\frac{\nabla\phi}{\phi}\nabla \frac{1}{\cos\theta}\notag ,
\end{align}
where  we used the Kato's inequality
\begin{equation*}
	|\nabla|H||^{2}\leq|\nabla^{\bot}H|^{2}.
\end{equation*}

Set $\frac{f}{f'}=g$. We choose $g$ such that for $x\in [1,\frac{1}{\delta}]$
\begin{eqnarray}\label{e6}
	\left\{\begin{array}{clcr} \varepsilon x/g>2,\\
		-4g'+8g/x-2=0,
	\end{array}\right.
\end{eqnarray}
where $\varepsilon$ is a fixed constant. Note that $x$ in this place represents the range of Lagrangian angle.
Let $g(x)=\frac{x}{2}c(x)$, then $c(x)$ need satisfy
\begin{eqnarray*}
	\left\{\begin{array}{clcr} 0<c(x)<\varepsilon,\\
		-xc'=1-c.
	\end{array}\right.
\end{eqnarray*}
We choose $c(x)=1-\eta x$ by solving the last equation. Then we have that
\begin{equation*}
	1-\frac{\eta}{\delta}\leq c(x)\leq1-\eta.
\end{equation*}
Therefore, the system (\ref{e6}) has a solution if and only if $1-\varepsilon<\eta<\delta$. Hence, our assumption that $\varepsilon>1-\delta$ implies that we can choose
\begin{equation*}
	\eta=\frac{1-\varepsilon+\delta}{2}.
\end{equation*}
Wth this range of $\eta$, we have
$$g=\frac{x}{2}(1-\eta x)$$ and
$$
f(x)=\frac{x^2}{(1-\eta x)^2},\quad x\in[1,
\frac{1}{\delta}].
$$
It is evident that for $x\in [1, \frac{1}{\delta}]$, $$\frac{1}{(1-\eta)^2}\leq
f(x)\leq \frac{1}{(\delta-\eta)^2}.
$$
On the other hand, by our assumption(\ref{varepsilon}), we know from (\ref{Hphi}) and the choice of $f$ that
\begin{equation}\label{Hphi2}
	\begin{split}
		\big(\Delta-\frac{\partial}{\partial t}\big) \phi
		&\geq 2\phi|A|^{2}\big(\frac{\varepsilon}{1-\frac{\eta}{\cos\theta}}-1\big)
		+\frac{\phi}{2}\big(\frac{|\nabla \phi|^{2}}{\phi^{2}}
		-2\frac{f'}{f}\nabla\frac{1}{\cos\theta}\cdot\frac{\nabla \phi}{\phi}\big) \\
		&\quad +2|H|^{2}f'\frac{\nabla\phi}{\phi}\cdot \nabla\frac{1}{\cos\theta} \\
		&\geq
		\frac{2(\varepsilon+\eta-1)}{1-\eta} \phi|A|^{2}+\frac{|\nabla \phi|^{2}}{2\phi}-
		\big(\phi\frac{f'}{f}\nabla\frac{1}{\cos\theta}-2|H|^{2}f'\nabla\frac{1}{\cos\theta}\big)\cdot
		\frac{\nabla\phi}{\phi}   \\
		& \geq \frac{2(\varepsilon+\delta-1)}{n(1+\varepsilon-\delta)} \phi|H|^{2} - \textbf{b} \cdot
		\frac{\nabla\phi}{\phi},
	\end{split}
\end{equation}
where
$\textbf{b}=\phi\frac{f'}{f}\nabla\frac{1}{\cos\theta}-2|H|^{2}f'\nabla\frac{1}{\cos\theta}$.

We claim that $\textbf{b}$ is bounded on $\Sigma_t$. Indeed, the definition of $f$ guarantees
that $f'$ and $\frac{f'}{f}$ are uniformly bounded by a constant
depending only on $\delta$. Therefore we only need to control
$\frac{\nabla
	\cos\theta}{\cos^{2}\theta}=\frac{-\sin\theta}{\cos^{2}\theta}\nabla
\theta$.
But by our assumption, $\cos\theta$ is bounded below by
$\delta$. Thus, (\ref{HJtheta}) gives the upper bound of $\frac{\nabla
	\cos\theta}{\cos^{2}\theta}$ because of the boundedness of
$|A|^{2}$.

Assume $|H|^{2}$ is not identically zero, we may assume $\phi(0,0)=1>0$.

Let $\psi(x)$ be a $C^{2}$ function on $[0,\infty)$ such that
\begin{eqnarray*}
	\psi(x)=\left\{\begin{array}{clcr}  1, &{\rm if} & 0\leq x\leq\frac{1}{2},\\
		0, &{\rm if}  & x\geq 1,
	\end{array}\right.
\end{eqnarray*}
$$ 0\leq\psi(x)\leq 1,~~~ \psi'(x)\leq 0,~~~ |\psi''(x)|\leq C ~~~{\rm and
}~~~ \frac{|\psi'(x)|^2}{\psi(x)}\leq C,
$$ 
where $C$ is an absolute constant.

We construct a cutoff function $h(X,t)=\psi(\frac{|X|^{2}}{R^{2}})\psi(\frac{t}{T})$ for any fixed positive number $T$.

Now we consider the function $F$ on $B_R(0)\times[s_{i},0]$ defined by
$$F(X,t)=h\phi(X,t).$$
Using the fact that $|\nabla X|^{2}=n$, a straightforward computation shows that
\begin{align}\label{e-h-1}
		(\Delta-\frac{\partial}{\partial t})h&=4\psi^{''}(\frac{|X|^{2}}{R^{2}})\psi(\frac{t}{T})\frac{\left\langle X,\nabla X\right\rangle^{2}}{R^{4}}
		+\psi(\frac{t}{T})\frac{2n\psi^{'}(\frac{|X|^{2}}{R^{2}})}{R^{2}}
		-\frac{1}{T}\psi^{'}(\frac{t}{T})\psi(\frac{|X|^{2}}{R^{2}})\nonumber\\
		&\geq-\frac{C}{R^{2}}-\frac{C}{T}
\end{align}
\begin{equation}\label{Hgnablag}
	\frac{|\nabla h|^2}{h}\leq\frac{C}{R^2}.
\end{equation}
Since $F\big|_{\partial
	(\overline{B}_R)\times[s_{i},0]}=0$, then $F$ must achieve its maximum at some $(Y,t)\in B_R(0)\times[s_{i},0]$.

By the maximum principle, at $(Y,t)$, we obtain that
\begin{equation*}
	\frac{\partial}{\partial t}F(Y,t)\geq0,
\end{equation*}
\begin{equation*}
	\nabla F(Y,t)=0,
\end{equation*}
\begin{equation*}
	\Delta F(Y,t)\leq0.
\end{equation*}
Hence,
\begin{equation*}
	\big(\Delta-\frac{\partial}{\partial t}\big)F=\big(\Delta-\frac{\partial}{\partial t}\big)h\phi\leq0,
\end{equation*}
\begin{equation}\label{e-h}
	\nabla h=-\frac{h}{\phi}\nabla\phi.
\end{equation}
Using (\ref{Hphi2}), (\ref{e-h-1}), (\ref{Hgnablag}) and (\ref{e-h}), we obtain
\begin{equation}\label{Hphi3}
	\begin{split}
		0&\geq\big(\Delta-\frac{\partial}{\partial t}\big)(h\phi) \\
		&=\phi\big(\Delta-\frac{\partial}{\partial t}\big)h+h\big(\Delta-\frac{\partial}{\partial t}\big)\phi+2\nabla h\cdot\nabla\phi \\
		&\geq  -\frac{C}{R^{2}}\phi-\frac{C}{T}\phi+\frac{2(\varepsilon+\delta-1)}{n(1+\varepsilon-\delta)}\phi|H|^{2}h - \textbf{b} \cdot
		\frac{\nabla\phi}{\phi}h+2\nabla h\cdot\big(-\frac{\phi}{h}\big)\nabla h  \\
		&=-\frac{C}{R^{2}}\phi-\frac{C}{T}\phi+\frac{2(\varepsilon+\delta-1)}{n(1+\varepsilon-\delta)}\phi|H|^{2}h + \textbf{b} \cdot\nabla h-2\frac{\phi}{h}|\nabla h|^{2} \\
		&\geq -\frac{C}{R^{2}}\phi-\frac{C}{T}\phi+\frac{2(\varepsilon+\delta-1)}{n(1+\varepsilon-\delta)}\phi|H|^{2}h-|\textbf{b}|\frac{C}{R}-2\frac{C}{R^{2}}\phi \\
		&\geq\frac{2(\varepsilon+\delta-1)}{n(1+\varepsilon-\delta)}|H|^{2}-\frac{C}{R}-\frac{C}{R^{2}}-\frac{C}{T},
	\end{split}
\end{equation}
where $C$ is constant which depend only on $\delta$ and the bound of $|A|^{2}$. 
By (\ref{Hphi3}), it implies that
\begin{equation*}
	\frac{2(\varepsilon+\delta-1)}{n(1+\varepsilon-\delta)}|H|^{2}(y,t)\leq\frac{C}{R}+\frac{C}{R^{2}}+\frac{C}{T},
\end{equation*}
which gives us the desired contradiction by letting $R\to+\infty$ and $T\to+\infty$.
This proves that the eternal solution satisfies
\begin{equation*}
	|H|^{2}\equiv0,
\end{equation*}
which implies that $\Sigma_t$ is a minimal submanifold.
Then by our assumption (\ref{varepsilon}), we know $|A|=0$. This finishes the proof of the theorem.
\hfill $\boxed{}$

\vspace{.2in}

\section{Rigidity theorems for Lagrangian translating soliton with zero Maslov class }

\vspace{.1in}

In this section, we will prove the rigidity theorems for Lagrangian translating soliton with zero Maslov class .

\vspace{.1in}

\noindent\textbf{Proof of Theorem \ref{thm1.3}:}
We prove it by contradiction.
Suppose there is a complete Lagrangian translating soliton $\Sigma$ with translating vector $T$ in $\mathbb{C}^{n}$ with $-D_{1}\leq\theta\leq D_{2}$, $|T|=1$ and 
\begin{equation}\label{infHa1}
	\inf_{\Sigma}|H|^{2}:=a>0,
\end{equation}
where $a\in(0,\frac{1}{2})$.
Set $u=\theta^{2}$, then by \eqref{Deltatheta}, we can easily see that
\begin{equation}\label{Deltatheta2}
	\Delta u=2|H|^{2}-\left\langle V, \nabla u\right\rangle,
\end{equation}
where $\Delta$ and $\nabla$ are the Laplacian and gradient operator on $\Sigma$ with respect to the induced metric, respectively. 

Let $\phi$ be any cutoff function on $\Sigma$ and $p$ be a positive number to be determined later. Multiplying both sides of \eqref{Deltatheta2} by $\phi^{2}u^{p}$ and integrating by parts yields
\begin{equation*}
	\begin{split}
		&2\int_{\Sigma}\phi^{2}u^{p}|H|^{2}d\mu-\int_{\Sigma}\phi^{2}u^{p}\left\langle V, \nabla u\right\rangle d\mu\\
		=&\int_{\Sigma}\phi^{2}u^{p}\Delta ud\mu
		=-p\int_{\Sigma}\phi^{2}u^{p-1}|\nabla u|^{2}d\mu-2\int_{\Sigma}\phi u^{p}\nabla\phi\cdot\nabla ud\mu
	\end{split}
\end{equation*}
Rearranging this equality and using Young’s inequality, we obtain (for any $\varepsilon>0$)
\begin{equation*}
	\begin{split}
		&2\int_{\Sigma}\phi^{2}u^{p}|H|^{2}d\mu+p\int_{\Sigma}\phi^{2}u^{p-1}|\nabla u|^{2}d\mu\\
		=&\int_{\Sigma}\phi^{2}u^{p}\left\langle V, \nabla u\right\rangle d\mu-2\int_{\Sigma}\phi u^{p}\nabla\phi\cdot\nabla ud\mu\\
		\leq&\int_{\Sigma}\phi^{2}u^{p}|V||\nabla u|d\mu+2\int_{\Sigma}\phi u^{p}|\nabla\phi||\nabla u|d\mu\\
		\leq&\varepsilon\int_{\Sigma}\phi^{2}u^{p+1}|V|^{2}d\mu+\frac{1}{4\varepsilon}\int_{\Sigma}\phi^{2}u^{p-1}|\nabla u|^{2}d\mu\\
		&+2\int_{\Sigma}\phi u^{p}|\nabla\phi||\nabla u|d\mu,
	\end{split}
\end{equation*}
which implies that
\begin{equation}\label{intphiuHV1}
	\begin{split}
		&\int_{\Sigma}\phi^{2}u^{p}(2|H|^{2}-\varepsilon u|V|^{2})d\mu+(p-\frac{1}{4\varepsilon})\int_{\Sigma}\phi^{2}u^{p-1}|\nabla u|^{2}d\mu\\
		\leq&2\int_{\Sigma}\phi u^{p}|\nabla\phi||\nabla u|d\mu.
	\end{split}
\end{equation}
From $1=|T|^{2}=|H|^{2}+|V|^{2}$, $-D_{1}\leq\theta \leq D_{2}$ and \eqref{infHa1}, we get that
\begin{equation*}
	2|H|^{2}-\varepsilon u|V|^{2}\geq 2a-\varepsilon D^{2}(1-a)\geq2a-\varepsilon D^{2},
\end{equation*}
where $D=\max\{D_{1},D_{2}\}$.
We first choose $\varepsilon=\frac{a}{D^{2}}$, then
\begin{equation*}
	2|H|^{2}-\varepsilon D^{2}|V|^{2}\geq a.
\end{equation*}
Then we obtain from \eqref{intphiuHV1} that
\begin{equation*}\label{intphiua}
	\begin{split}
		&a\int_{\Sigma}\phi^{2}u^{p}d\mu+(p-\frac{D^{2}}{4a})\int_{\Sigma}\phi^{2}u^{p-1}|\nabla u|^{2}d\mu\\
		\leq&2\int_{\Sigma}\phi u^{p}|\nabla\phi||\nabla u|d\mu.
	\end{split}
\end{equation*}
Next, we choose $p=\frac{D^{2}}{4a}+1:=P_{0}$, then
\begin{equation}\label{intphiub}
	\begin{split}
		a\int_{\Sigma}\phi^{2}u^{P_{0}}d\mu+\int_{\Sigma}\phi^{2}u^{P_{0}-1}|\nabla u|^{2}d\mu
		\leq2\int_{\Sigma}\phi u^{P_{0}}|\nabla\phi||\nabla u|d\mu.
	\end{split}
\end{equation}

Now we take a smooth function $\phi(x):\mathbb{R}\to\mathbb{R}$ supported on $[-R,R]$ which has the property:
$\phi=1$ on $[-\frac{R}{2},\frac{R}{2}]$ and there is a constant $C_{b}>0$ such that $|\nabla \phi|\leq C_{b}\phi^{b}\frac{1}{R}$ for $0<b<1$.

Recall that we assume the Lagrangian translating soliton to have Euclidean area growth.
By using Young’s inequality and the property of $\phi$, we can estimate $\int_{\Sigma}\phi u^{P_{0}}|\nabla\phi||\nabla u|d\mu$ as follows:
\begin{equation}\label{intphiud}
	\begin{split}
		\int_{\Sigma}\phi u^{P_{0}}|\nabla\phi||\nabla u|d\mu&
		\leq\frac{C_{b}}{R}\int_{\Sigma}\phi^{1+b}u^{P_{0}}|\nabla u|d\mu
		\\&\leq2D\frac{C_{b}}{R}\int_{\Sigma}\phi^{1+b}u^{\frac{P_{0}(1+b)}{2}}u^{\frac{P_{0}(1-b)}{2}}d\mu\\
		&\leq\frac{a}{4}\int_{\Sigma}\phi^{2}u^{P_{0}}d\mu+\frac{C}{R^{\frac{2}{1-b}}}\int_{\Sigma\cap supp\phi}u^{P_{0}}d\mu\\
		&\leq\frac{a}{4}\int_{\Sigma}\phi^{2}u^{P_{0}}d\mu+CD^{2P_{0}}R^{n-\frac{2}{1-b}},
	\end{split}
\end{equation}
where we used the fact 
$|\nabla u|=|\nabla\theta^{2}|=2|\theta||\nabla\theta|\leq2D.$
Finally we choose $b>0$ such that $n-\frac{2}{1-b}<0$, i.e., $b>1-\frac{2}{n}$.
Using \eqref{intphiud} in \eqref{intphiub}, we have
\begin{equation*}\label{intphiuf}
	\begin{split}
		\frac{a}{2}\int_{\Sigma}\phi^{2}d\mu \leq C(a,D,n)R^{n-\frac{2}{1-b}} .
	\end{split}
\end{equation*}
This gives the desired contradiction as $R\to +\infty$.
\hfill $\boxed{}$

\vspace{.1in}

\noindent\textbf{Proof of Theorem \ref{thm1.4}:}
Let $\phi$ be any cutoff function and $p$ be a positive number to be determined later. Multiplying both sides of \eqref{Deltatheta2} by $\phi^{2}ue^{\left\langle T,{\bf{x}}\right\rangle}$ and integrating by parts yields
\begin{equation*}
	\begin{split}
		2\int_{\Sigma}\phi^{2}u|H|^{2}e^{\left\langle T,{\bf{x}}\right\rangle}d\mu
		=&\int_{\Sigma}\phi^{2}u(\Delta u+\left\langle {V},\nabla u \right\rangle)e^{\left\langle T,{\bf{x}}\right\rangle}d\mu\\
		=&\int_{\Sigma}\phi^{2}u\text{div}_{\Sigma}(e^{\left\langle T,{\bf{x}}\right\rangle}\nabla u)d\mu\\
		=&-\int_{\Sigma}\left\langle \nabla(\phi^{2}u),\nabla u \right\rangle e^{\left\langle T,{\bf{x}}\right\rangle}d\mu\\
		=&-\int_{\Sigma}\phi^{2}|\nabla u|^{2} e^{\left\langle T,{\bf{x}}\right\rangle}d\mu\\
		&-2\int_{\Sigma}\phi u\left\langle\nabla\phi ,\nabla u \right\rangle e^{\left\langle T,{\bf{x}}\right\rangle}d\mu
	\end{split}
\end{equation*}
which implies that
\begin{equation}\label{intphiue1}
	\begin{split}
		\int_{\Sigma}\phi^{2}(2u|H|^{2}+|\nabla u|^{2})e^{\left\langle T,{\bf{x}}\right\rangle}d\mu
		\leq2\int_{\Sigma}\phi u|\nabla\phi||\nabla u| e^{\left\langle T,{\bf{x}}\right\rangle }d\mu.
	\end{split}
\end{equation}
Since $-D_{1}\leq\theta\leq D_{2}$, we see that $\theta\leq D:=\max\{D_{1},D_{2}\}$. 
Furthermore, from \eqref{nablatheta=H}, we have
\begin{equation}
	|\nabla u|=|\nabla\theta^{2}|=2|\theta||\nabla\theta|\leq2D|H|.
\end{equation}
Therefore, by \eqref{intphiue1}, we have
\begin{equation}\label{intphiuD}
	\begin{split}
		\int_{\Sigma}\phi^{2}(2u|H|^{2}+|\nabla u|^{2})e^{\left\langle T,{\bf{x}}\right\rangle}d\mu
		&\leq C(D)\int_{\Sigma}\phi|\nabla\phi||H| e^{\left\langle T,{\bf{x}}\right\rangle }d\mu\\
		&=C(D)\int_{\Sigma}\phi|\nabla\phi||H|d\widetilde{\mu},
	\end{split}
\end{equation}
where recall $d\widetilde{\mu}=e^{\left\langle T,{\bf{x}}\right\rangle}d\mu$.
Now for any fixed $R>0$, we take the cutoff function $\phi=\phi_{R}$ such that $\phi\in C^{\infty}_{0}(B(2R))$,
$\phi=1$ on $B(R)$, and $|\nabla\phi|\leq|D\phi|\leq\frac{C}{R}$.
Here, $B(R)$ is the ball of radius in $R^{n}$, $D\phi$ is the gradient of $\phi$ with respect to the Euclidean metric in ${\mathbb R}^{n}$, and $C$ is an absolute constant. 
Taking $\phi=\phi_{R}$ in \eqref{intphiuD} yields
\begin{equation}\label{intphiuDW}
	\begin{split}
		\int_{\Sigma}\phi^{2}(2u|H|^{2}+|\nabla u|^{2})e^{\left\langle T,{\bf{x}}\right\rangle}d\mu
		&\leq C(D)\int_{\Sigma}\phi|\nabla\phi||H|d\widetilde{\mu}\\
		&\leq \frac{C(D)}{R}\int_{\Sigma}|H|d\widetilde{\mu}.
	\end{split}
\end{equation}
By our assumption, $\int_{\Sigma}|H|d\widetilde{\mu}<\infty$.
Letting $R\to+\infty$ in \eqref{intphiuDW}, we finally obtain that $H\equiv0$ on $\Sigma$. 
This means that $\Sigma$ is a minimal submanifold.

When $n=2$, by \eqref{AHV}, we see that $A\equiv0$ on $\Sigma$.
Thus it must be a flat plane.
\hfill $\boxed{}$

\vspace{.1in}

\noindent\textbf{Proof of Theorem \ref{thm1.5}:}

Let $\phi$ be any cutoff function and $p$ be a positive number to be determined later. Multiplying both sides of \eqref{Deltatheta2} by $\phi^{2}ue^{\left\langle T,{\bf{x}}\right\rangle}$ and integrating by parts yields
\begin{equation*}
	\begin{split}
		2\int_{\Sigma}\phi^{2}u|H|^{2}e^{\left\langle T,{\bf{x}}\right\rangle}d\mu
		=&\int_{\Sigma}\phi^{2}u(\Delta u+\left\langle V,\nabla u \right\rangle)e^{\left\langle T,{\bf{x}}\right\rangle}d\mu\\
		=&\int_{\Sigma}\phi^{2}u\text{div}_{\Sigma}(e^{\left\langle T,{\bf{x}}\right\rangle}\nabla u)d\mu\\
		=&-\int_{\Sigma}\left\langle \nabla(\phi^{2}u),\nabla u \right\rangle e^{\left\langle T,{\bf{x}}\right\rangle}d\mu\\
		=&-\int_{\Sigma}\phi^{2}|\nabla u|^{2} e^{\left\langle T,{\bf{x}}\right\rangle}d\mu\\
		 &-2\int_{\Sigma}\phi u\left\langle\nabla\phi ,\nabla u \right\rangle e^{\left\langle T,{\bf{x}}\right\rangle}d\mu
	\end{split}
\end{equation*}
which implies that
\begin{equation}\label{intphiue}
	\begin{split}
		\int_{\Sigma}\phi^{2}(2u|H|^{2}+|\nabla u|^{2})e^{\left\langle T,{\bf{x}}\right\rangle}d\mu
		\leq2\int_{\Sigma}\phi u|\nabla\phi||\nabla u| e^{\left\langle T,{\bf{x}}\right\rangle }d\mu.
	\end{split}
\end{equation}

Now we take a smooth function $\phi(x):\mathbb{R}\to\mathbb{R}$ supported on $[-R,R]$ which has the property:
$\phi=1$ on $[-\frac{R}{2},\frac{R}{2}]$ and there is a constant $C_{b}>0$ such that $|\nabla \phi|\leq C_{b}\phi^{b}\frac{1}{R}$ for $0<b<1$.

By using Young’s inequality, the bound of $u=\theta^{2}$ and the property of $\phi$, we can estimate \eqref{intphiue} as follows:
\begin{equation*}
	\begin{split}
		\int_{\Sigma}\phi^{2}(2u|H|^{2}+|\nabla u|^{2})e^{\left\langle T,{\bf{x}}\right\rangle}d\mu
		\leq&C\int_{\Sigma}\phi^{1+b}u|\nabla u|\frac{1}{R}e^{\left\langle T,{\bf{x}}\right\rangle }d\mu\\
		\leq&\widetilde{C}\frac{1}{R}\int_{\Sigma}\phi^{1+b}|\nabla u|e^{\left\langle T,{\bf{x}}\right\rangle }d\mu\\
		\leq&\int_{\Sigma}\phi^{2}|\nabla u|^{2}e^{\left\langle T,{\bf{x}}\right\rangle }d\mu
		+\int_{\Sigma\cap supp\phi}\frac{C}{R^{\frac{2}{1-b}}}e^{\left\langle T,{\bf{x}}\right\rangle }d\mu.
	\end{split}
\end{equation*}
Finally we choose $b>0$ such that $n-\frac{2}{1-b}<0$, i.e., $b>1-\frac{2}{n}$.
Then
\begin{equation*}\label{intphiuemu}
	\begin{split}
		\int_{\Sigma}\phi^{2}u|H|^{2}e^{\left\langle T,{\bf{x}}\right\rangle}d\mu
		\leq\frac{C}{R^{\frac{2}{1-b}}}\widetilde{\mu}(\Sigma\cap B_{2R}(0)).
	\end{split}
\end{equation*}
Letting $R\to+\infty$ in \eqref{intphiuemu}, we obtain that $|H|\equiv0$ on $\Sigma$.
This means that $\Sigma$ is a minimal submanifold.

When $n=2$, by \eqref{AHV}, we see that $A\equiv0$ on $\Sigma$.
Thus it must be a flat plane.
\hfill $\boxed{}$

\vspace{.3in}

\textbf{Acknowlegement.}  The second author was supported by the National Natural Science Foundation of China (Grant No. 12271069) and the third author was supported by the National Natural Science Foundation of China (Grant No. 12071352, 12271039). The authors would like to thank Professor Hongbing Qiu for stimulating discussions on this subject. Many thanks to the referee for pointing out a lot of typos and mathematical suggestions which made this paper much more readable.

\vspace{.2in}

\noindent\textbf{Author Contributions} All authors have contributed equally.

\vspace{.2in}

\noindent\textbf{Data availability} No data sets were generated or analysed during the current study.
\vspace{.2in}

\noindent{\Large\textbf{Declarations}}

\vspace{.2in}

\noindent\textbf{Conflict of interest} Not applicable.

\end{document}